# BAYESIAN VARIABLE SELECTION AND DATA INTEGRATION FOR BIOLOGICAL REGULATORY NETWORKS[1]

By Shane T. Jensen,[2] Guang Chen[3] and Christian J. Stoeckert, Jr.[3]

*University of Pennsylvania*

A substantial focus of research in molecular biology are *gene regulatory networks*: the set of transcription factors and target genes which control the involvement of different biological processes in living cells. Previous statistical approaches for identifying gene regulatory networks have used gene expression data, ChIP binding data or promoter sequence data, but each of these resources provides only partial information. We present a Bayesian hierarchical model that integrates all three data types in a principled variable selection framework. The gene expression data are modeled as a function of the unknown gene regulatory network which has an informed prior distribution based upon both ChIP binding and promoter sequence data. We also present a variable weighting methodology for the principled balancing of multiple sources of prior information. We apply our procedure to the discovery of gene regulatory relationships in *Saccharomyces cerevisiae* (Yeast) for which we can use several external sources of information to validate our results. Our inferred relationships show greater biological relevance on the external validation measures than previous data integration methods. Our model also estimates synergistic and antagonistic interactions between transcription factors, many of which are validated by previous studies. We also evaluate the results from our procedure for the weighting for multiple sources of prior information. Finally, we discuss our methodology in the context of previous approaches to data integration and Bayesian variable selection.

Received April 2007; revised August 2007.

[1]The first two authors contributed equally to this work.

[2]Supported in part by grant from the University of Pennsylvania Research Foundation.

[3]Supported in part by NIH Grant U01-DK56947.

Supplementary material available at http://imstat.org/aoas/supplements

**1. Introduction and motivation.**    The development and function of living cells is, to a large extent, dictated by a carefully choreographed system of gene expression. Gene expression is controlled in part by *transcription factors* (TFs), a class of proteins which bind to DNA leading to an increase or decrease in transcription of target genes. The collection of transcription factors and their targets (genes that they control) is called a regulatory network. In this work, we develop a model for understanding the transcriptional regulatory networks that specify and maintain cellular function.

The computational approaches that are used to identify regulatory networks have typically used information from three different sources:

1. *Gene expression data.* Microarray chips are used to measure the levels of mRNA produced for each gene in a cell, which is usually referred to as the amount of gene expression. Since mRNA is a precursor to the protein product of each gene, expression levels are used as proxy for the amount of protein produced. Genes which show similar levels of expression in different conditions are believed to be co-regulated.
2. *ChIP binding data.* Chromatin Immunoprecipitation technology uses antibodies to isolate sequences that are directly bound by a specific transcription factor. Microarray chips are then used to chart these sequences within the genome in order to determine potential locations for binding of that particular transcription factor.
3. *Promoter sequence data.* The different binding sites (located near different target genes) of the same transcription factor show a significant sequence conservation, but substantial variability is also present. The conserved appearance of the transcription factor binding sites is summarized by a position-specific weight matrix (PWM) which can be used to search near to potential target genes for the predicted binding sites of a transcription factor. The strength of the signal in these PWMs varies substantially between transcription factors.

Although each of these resources are extremely useful, their power is inherently limited by the fact that each type of data provides only partial information: expression data provide only indirect evidence of regulation, promoter sequence data provide only potential binding sites which may not be bound by TFs and ChIP binding data provide only physical binding locations which may not be functional in terms of controlling gene expression. We develop a Bayesian hierarchical model for combining our three available sources of information: gene expression data, ChIP binding data and promoter sequence data. This is accomplished by extending previous linear models for gene expression data [Bussemaker, Li and Siggia (2001), Gao, Foat and Bussemaker (2004)] into a variable selection framework. There has been substantial research into Bayesian approaches to variable selection, though as mentioned by George (2000), most previous meth-ods have focused on minimization of prior dependence. In contrast, our approach takes advantage of two additional data sources, ChIP binding and promoter element data, to generate informed prior distributions for our variable selection model. We also develop a variable weighting methodology for balancing these two sources of prior information.

There has also been substantial previous research into the integration of biological data sources for the discovery of regulatory networks. Bussemaker, Li and Siggia (2001) developed a linear model to reflect the correlation between expression patterns and cis-regulatory motif abundance, with the inherent drawback that any synergistic effects from transcription factor interaction were not taken into account. Tadesse, Vannucci and Lio (2004) posited a similar linear model between expression patterns and cis-regulatory motif abundance, but used Bayesian variable selection instead of the stepwise regression procedure of Bussemaker, Li and Siggia (2001). Gao, Foat and Bussemaker (2004) presented an integrated linear model, MA-Networker, for combining expression and ChIP binding data, but their procedure required a stringent binding p-value threshold. Banerjee and Zhang (2003) used thresholded ChIP binding data and gene expression data to identify cooperativity among TFs. Our model also allows us to estimate synergistic and antagonistic interactions between transcription factors. Xing and Laan (2005) developed a multiple linear regression model selected by a loss-based V-fold cross-validation selector, but the method relies on knowledge of known TF sites and the number of TFs with known consensus binding sites is small and their functional coverage is somewhat limited. Based on the assumption that the expression levels of regulated genes depend on the expression levels of regulators, Segal (2001, 2003) constructed a probabilistic model which used binding motif features and expression data to identify modules of co-regulated genes and their regulators. This probabilistic model reflected nonlinear properties, but required prior clustering of the expression data. The GRAM algorithm combining ChIP binding and expression data was developed by Bar-Joseph et al. (2003) to discover regulatory networks in *Saccharomyces cerevisiae*, but their technique is heuristic with arbitrary parameter thresholds and little systematic modeling. Another probabilistic-based approach is taken by Lemmens et al. (2006) in their ReMoDiscovery method.

Our full probabilistic model does not rely on pre-clustering of expression data and reduces dependence on arbitrary parameter cutoffs. As mentioned by Kloster, Tang and Wingreen (2005), many current methods rely on the basic assumption that each gene can only belong to a single cluster. Our framework permits genes to belong to multiple regulatory clusters, which allows us to model multiple biological pathways simultaneously. Other recent efforts [Liao et al. (2003), Yang (2005), Boulesteix and Strimmer (2005)] have used a "network component analysis" approach to find regulatory modules using expression data. In these investigations, ChIP binding data are used to form a connectivity network between genes and TFs, which is assumed to be known without error. In contrast, our model allows for the inherent uncertainty in ChIP experiments, which allows for a more direct integration of ChIP binding and gene expression, but uses TF gene expression as our measure of TF activity. Sabatti (2005) also use the "network component analysis" approach to model gene expression using a prior distribution based on promotor binding sites, which did allow for some uncertainty in the sequence information, but did not include any ChIP binding data. We will revisit the distinction between our approach and "network component analysis" in our discussion.

In Section 2 we outline our Bayesian variable selection model for integrating multiple data sources and discuss implementation using a Gibbs sampling algorithm [Geman (1984)]. In Section 3 we describe an application of our methodology to *Saccharomyces cerevisiae* (Yeast) where gene expression, ChIP binding data and promoter sequence data are available for many transcription factors. We validate our results in yeast using external information from the biological literature and compare to several alternative methods. In a related paper [Chen, Jensen and Stoeckert (2007)] we present a reduced application of our model to Yeast, as well as applications in the higher organism *Mus musculus* (mouse). We explore several interesting model consequences and sensitivities in Section 4. Finally, in Section 5 we discuss our model in the context of previous methods for regulatory network elucidation, as well as previous approaches to variable selection.

## 2. Bayesian model and implementation.

The primary goal of our statistical model is to infer probable gene–TF relationships through the integration of available biological data. Mathematically, we formulate these relationships as unknown indicator variables $C_{ij} = 1$ if gene $i$ is regulated by TF $j$ or 0 otherwise. Our inference for these regulation indicators $C_{ij}$ is a variable selection process that determines which subset of the many possible gene–TF relationships are biologically important and allows us to construct an inferred regulatory network. This network can be visually represented as a graph where nodes are genes and TFs, and each $C_{ij}$ variable determines whether or not there should be a directed edge connecting the node for TF $j$ with the node for gene $i$. Collectively, the matrix $\boldsymbol{C}$ of these indicator variables also gives us regulatory clusters (also called regulatory modules) for each TF, since all genes $i$ where $C_{ij} = 1$ are estimated to be in a cluster together regulated by TF $j$. An important aspect of our flexible framework is that we are explicitly allowing genes to belong to multiple clusters controlled by different transcription factors (i.e., $C_{ij} = 1$ and $C_{ij'} = 1$ for $j \neq j'$). In order to infer likely values for our indicator variables $\boldsymbol{C}$, our model incorporates up to three general classes of biological information: gene expression data, ChIP binding data and sequence-level promoter data.

We denote our gene expression data as $g_{it}$, the expression of gene $i$ ($i = 1, \ldots, N$) in experiment $t$ ($t = 1, \ldots, T$). The set of $T$ experiments can be from different tissues, time-course experiments, different gene-knockout experiments, or any combination thereof. Within these expression data, we give special focus to the expression of genes that produce known transcription factor proteins. For our $J$ known transcription factors, we denote $f_{jt}$ as the expression of TF $j$ ($j = 1, \ldots, J$) in experiment $t$. The TF expression levels $f_{jt}$ are derived from our gene expression data by simply identifying the gene that encodes each transcription factor $j$, and using the expression level of that gene as our TF expression levels. We will use the expression $f_{jt}$ for the gene that produces TF $j$ as a proxy for the amount of activity of TF $j$. In addition to expression data, we have available Chromatin Immunoprecipitation (ChIP) experiments which give information on the physical binding location of specific transcription factors. We use $b_{ij}$ to denote the probability that transcription factor $j$ physically binds in close proximity to gene $i$, from a ChIP binding experiment for transcription factor $j$. Finally, we have available sequence-level information in the form of known or putative promoter binding sites for specific transcription factors located in the upstream regions of target genes. We denote $m_{ij}$ as the probability that transcription factor $j$ has a promoter binding site in the regulatory region of gene $i$. These binding sites could be experimentally verified or predicted by scanning upstream sequences for similarity to an established position-specific weight matrix (PWM) for a particular transcription factor. We outline our model in the most general case where all three of these data types are present, but we will also discuss the ramifications on our procedure when only subsets of these data types are available. Our different data sources are summarized in Table 1.

The first level of our probabilistic model incorporates our gene expression data by specifying the observed gene expression $g_{it}$ as a linear function of TF expression, $f_{jt}$,

$$(1) \qquad g_{it} = \alpha_i + \sum_{j=1}^{J} \beta_j C_{ij} f_{jt} + \epsilon_{it}, \qquad \epsilon_{it} \sim \text{Normal}(0, \sigma^2).$$

TABLE 1
*Notation for available data sources*

| Notation | Data type |
| --- | --- |
| $g_{it}$ = expression of gene $i$ in experiment $t$ | Gene expression |
| $f_{jt}$ = expression of TF $j$ in experiment $t$ | TF expression |
| $b_{ij}$ = probability that TF $j$ binds near to gene $i$ | ChIP binding |
| $m_{ij}$ = probability that gene $i$ has promoter element for TF $j$ | Promoter sequence |

In equation (1) we see that our regulation indicators $C_{ij}$ act as variable selection parameters: only TFs $j$ where $C_{ij} = 1$ are allowed to influence the expression of gene $i$. The parameter $\beta_j$ is the linear effect of TF $j$ on gene expression, whereas $\alpha_i$ can be interpreted as the baseline expression for gene $i$ in absence of regulation by known transcription factors (i.e., $C_{ij} = 0$ for all TFs $j$). Bussemaker, Li and Siggia (2001) also used a linear model for expression data, except that their approach did not use TF expression $f_{jt}$ as a proxy for TF activity, but rather used sequence elements as their proxy for TF activity. We prefer the use of TF expression as our proxy for TF activity since our TF expression levels $f_{jt}$ are specific to each experiment $t$ in the same way as our gene expression levels $g_{it}$. Sequence information is not experiment or condition-specific and so is less useful as a proxy for TF activity. However, we do make use of sequence elements in our prior distribution (3) for the global regulation indicators $C_{ij}$.

Our simple linear model, as stated in (1), is limited by not allowing for combinatorial relationships between TFs. Each TF $j$ has a single effect ($\beta_j$) on the expression of gene $i$, which does not take into account the biological reality that expression is often the result of synergistic or antagonistic action of multiple TFs binding simultaneously. We acknowledge these combinatorial relationships by expanding our linear model to include interaction terms:

$$(2) \qquad g_{it} = \alpha_i + \sum_{j=1}^{J} \beta_j C_{ij} f_{jt} + \sum_{j \neq k} \gamma_{jk} C_{ij} C_{ik} f_{jt} f_{kt} + \epsilon_{it},$$

$$\epsilon_{it} \sim \text{Normal}(0, \sigma^2),$$

where we now have additional coefficients $\gamma_{jk}$ that can be interpreted as the synergistic (or antagonistic) effect of both TFs $j$ and $k$ binding together to the same upstream region (in addition to the effects of TF $j$ or $k$ binding in isolation). Note that our regulation indicators $C_{ij}$ again act as variable selectors for both the linear and interaction terms in equation (2). Of course, higher-order interactions or nonlinear functions could also be considered in our framework. However, this additional model complexity would increase the parameter space and computation burden of the model dramatically. We believe that our extended model with TF interactions (2) achieves an appropriate balance between the computational cost of model fitting and the flexibility to adequately model TF-gene expression relationships.

Despite the intuitive appeal of positing linear models [Bussemaker, Li and Siggia (2001), Tadesse, Vannucci and Lio (2004), Gao, Foat and Bussemaker (2004)] as a variable selection problem, the implementation of our variable selection model is quite complex in practice, with a large number of both genes $i$ (e.g., 6026 in our yeast application) and TFs $j$ (e.g., 39 in our yeast application). We address this complexity by using our additional data types to construct informed prior distributions for each regulation indicator $C_{ij}$. We have $b_{ij}$, the probability that TF $j$ physically binds in proximity to gene $i$ in a ChIP-binding experiment, and $m_{ij}$, the probability of a binding site for TF $j$ in the upstream region of gene $i$. The second component of our model incorporates both $b_{ij}$ and $m_{ij}$ into a combined prior distribution for our unknown regulation indicators $C_{ij}$:

$$(3) \qquad p(C_{ij}|m_{ij}, b_{ij}, w_j) \propto [b_{ij}^{C_{ij}}(1-b_{ij})^{1-C_{ij}}]^{w_j} \cdot [m_{ij}^{C_{ij}}(1-m_{ij})^{1-C_{ij}}]^{1-w_j}.$$

The variable $w_j$ is the relative weight of the prior ChIP-binding information $b_{ij}$ versus the TF binding site information $m_{ij}$. The weights $\boldsymbol{w} = (w_1, \ldots, w_J)$ are TF-specific but not gene-specific, and are designed to reflect potential global differences in quality between the binding data and promoter sequence data for TF $j$. However, since this relative quality is not necessarily known a priori, we will treat each weight $w_j$ as an unknown variable. Clearly, if only ChIP binding data for TF $j$ are available, then $w_j = 1$ and equation (3) reduces to a function of $b_{ij}$ only, whereas if only promoter sequence data for TF $j$ are available, then $w_j = 0$ and equation (3) reduces to a function of $m_{ij}$ only. In cases where both data types are available, our model will estimate the weight $w_j$ so that our prior distribution moves toward our likelihood based on expression data, thereby creating an appropriate balance between the two sources of prior information.

The Bayesian approach gives us a principled framework for connecting these model components into a single posterior distribution for all unknown parameters:

$$p(\boldsymbol{C}, \boldsymbol{w}, \boldsymbol{\Theta}|\boldsymbol{g}, \boldsymbol{f}, \boldsymbol{m}, \boldsymbol{b}) \propto p(\boldsymbol{g}|\boldsymbol{f}, \boldsymbol{C}, \boldsymbol{\Theta}) \cdot p(\boldsymbol{C}|\boldsymbol{m}, \boldsymbol{b}, \boldsymbol{w}) \cdot p(\boldsymbol{\Theta}, \boldsymbol{w}),$$

where $\boldsymbol{\Theta}$ denotes the collection of linear model parameters, that is, $\boldsymbol{\Theta} = (\boldsymbol{\alpha}, \boldsymbol{\beta}, \boldsymbol{\gamma}, \sigma^2)$. The term $p(\boldsymbol{g}|\boldsymbol{f}, \boldsymbol{C}, \boldsymbol{\Theta})$ represents our first model level with expression data $\boldsymbol{g} = (g_{it})$ and $\boldsymbol{f} = (f_{jt})$ and $p(\boldsymbol{C}|\boldsymbol{m}, \boldsymbol{b}, \boldsymbol{w})$ represent our second model level with ChIP binding data $\boldsymbol{b} = (b_{ij})$ and promoter sequence data $\boldsymbol{m} = (m_{ij})$. All that remains is the specification $p(\boldsymbol{\Theta}, \boldsymbol{w})$, the prior distributions for our TF-specific prior weights $\boldsymbol{w} = (w_1, \ldots, w_J)$ and our linear model parameters $\boldsymbol{\Theta}$:

(a) baseline gene $i$ expression: $\alpha_i \sim \text{Normal}(0, \tau_\alpha^2)$,
(b) TF linear effects: $\beta_j \sim \text{Normal}(0, \tau_\beta^2)$,
(c) TF interaction effects: $\gamma_{jk} \sim \text{Normal}(0, \tau_\gamma^2)$,
(d) residual gene expression variance: $\sigma^2 \sim \text{Inv-}\chi_\nu^2$,
(e) prior distribution weights: $w_j \sim \text{Uniform}(0, 1)$.

In Section 4.2 we discuss choices of these hyper-parameters $\boldsymbol{\tau}$ and $\nu$ that are noninfluential on our posterior inference. We estimate the joint posterior distribution of all unknown parameters by Markov chain Monte Carlo simulation. Specifically, we use Gibbs sampling [Geman (1984)], where we iteratively sample values of one set of parameters given all other parameters:

1. Sampling $\boldsymbol{\Theta}$ given $\boldsymbol{C}, \boldsymbol{w}$ and data $\boldsymbol{g}, \boldsymbol{f}, \boldsymbol{b}, \boldsymbol{m}$.
2. Sampling $\boldsymbol{C}$ given $\boldsymbol{w}, \boldsymbol{\Theta}$ and data $\boldsymbol{g}, \boldsymbol{f}, \boldsymbol{b}, \boldsymbol{m}$.
3. Sampling $\boldsymbol{w}$ given $\boldsymbol{C}, \boldsymbol{\Theta}$ and data $\boldsymbol{g}, \boldsymbol{f}, \boldsymbol{b}, \boldsymbol{m}$.

The details of our Gibbs sampling implementation are given in the Appendix. Software for our procedure is available for download at http://www.cbil.upenn.edu/COGRIM/.

**3. Application to the yeast regulatory network.** We applied our model to extensive available data for the simple organism *Saccharomyces cerevisiae* (budding yeast). We used 314 gene expression experiments, each involving 6026 yeast genes. A detailed reference list for our expression data sources and description of some preliminary data cleaning and manipulation is given in the supplemental materials. We also have both ChIP binding data [Lee et al. (2002)] and promoter element data [Matys et al. (2003)] for 39 yeast transcription factors. Thus, the dimension of our observed expression data $\boldsymbol{g} = (g_{it})$ is 6026 genes × 314 conditions, while $\boldsymbol{f} = (f_{jt})$ is 39 TFs × 314 conditions. The dimensions of our observed binding data $\boldsymbol{b} = (b_{ij})$ and $\boldsymbol{m} = (m_{ij})$ are both 6026 genes × 39 TFs. Our supplemental materials also contain a detailed evaluation of the convergence of our Gibbs sampling algorithm.

In Section 3.1 below we examine our posterior results for the regulation indicators $\boldsymbol{C}$ in our model, which are the primary goal of our investigation. We use available external information from the biological literature to confirm our inference and compare to previous methods. In Section 3.2 we present additional results for the interaction between TFs in our yeast application where the regulatory actions of many TFs are being modeled simultaneously. Finally, in Section 3.3 we examine posterior inference for our weighting parameters between the two different sources of available prior information for each Yeast transcription factor. A reduced form of our model is applied to Yeast and two transcription factors in the higher organism *Mus musculus* (mouse) in a related paper [Chen, Jensen and Stoeckert (2007)].

3.1. *Inference for regulation indicators* $\boldsymbol{C}$. The samples of each indicator variable $C_{ij}$ from our Gibbs sampling algorithm were used to estimate the posterior probability $P(C_{ij} = 1)$ for each possible gene $i$ and TF $j$ relationship. We considered any $(i, j)$ combination with posterior probability $P(C_{ij} = 1)$ higher than 0.5 as an inferred gene–TF relationship, and we then call the $i$th gene a *target gene* of transcription factor $j$. For our yeast application, we focus on the inferred target genes for 39 transcription factors where external validation measures of biological relevance are available. We also give a visual representation of the regulatory network in our supplemental materials. We use two validation measures, TF knockout data and MIPS functional categories, to compare the inferred gene sets from our full model involving all three data sources to reduced forms of our model that only involve subsets of our data sources. With these same two validation measures, we also compare our model inference to the inferred target genes produced by several previous integration methods: GRAM [Bar-Joseph et al. (2003)], ReMoDiscovery [Lemmens et al. (2006)] and MA-Networker [Gao, Foat and Bussemaker (2004)]. Finally, we compare our model results to target gene sets constructed based on heuristic thresholds of single data sources used in isolation. Following the recommendation of Lee et al. (2002), we use thresholded ChIP-binding data alone by classifying any genes with binding p-values less than 0.001 as gene targets. We also use thresholded expression data alone by calculating the pairwise correlation between gene expression $\boldsymbol{g}_i$ and the expression $\boldsymbol{f}_j$ of TF $j$, and classifying the most correlated 1% of genes as targets. This 1% threshold gave the best performance among several different thresholds that we considered.

Our most reliable validation uses the results of TF *knockout* experiments from the Rosetta Yeast Compendium [Hughes et al. (2000)] for four yeast TFs: Yap1, Swi4, Swi5 and Gcn4. Knockout experiments are considered a gold standard for the regulatory activity of individual transcription factors. In each of these experiments, a knockout strain of yeast was created with a specific TF removed from the genome. Microarray chips are then used to quantify the *knockout response* for each gene: the change in expression for each gene between the knockout and wild-type strains. Genes that are targets of the knocked-out TF should show greater knockout response between the wild-type strain and the knock-out strain.

Within each TF knockout experiment, we calculated a t-statistic for the knockout response for genes inferred to be targets by each method, which are shown in Figure 1. Methods with larger t-statistic values in Figure 1 show a greater knockout response within their inferred target genes, which supports the biological relevance of that method. For each TF experiment, our model using expression data only ("Exp") is clearly inferior to our model with multiple data sources ("All 3" and "ExpChIP"). Our model based on all three data sources ("All 3") shows similar performance to our model without promoter sequence data ("Exp+ChIP"), which suggests that this third data source is not contributing substantially to inference. We will revisit this issue when we examine our variable weight inference in Section 3.3. The inferred target genes from our integrated models ("All 3" and "ExpChIP") show uniformly superior performance across the four experiments, suggesting that our full probabilistic model is capturing more signal than previous integrated methods (MA-Networker, GRAM and ReMoDiscovery).

The inferred target genes based on thresholded single data shows considerably worse performance, which demonstrates that an integrated approach based on multiple data sources leads to superior inference for regulatory networks.

As a second validation, we used the MIPS database [Mewes (2002)] to assign a functional category to each gene in our yeast application. For each of our 39 TFs, we looked for over-represented functional categories within each set of inferred target genes. A set of putative gene targets that share similar gene functions are likely to be involved in the same biological pathway, which validates the inference that they are regulated by a common transcription factor. Any functional category with a p-value of less than 0.001 (p-value calculated using the hypergeometric distribution) was considered to be significantly over-represented. The proportion of inferred target genes that shared over-represented functional categories (averaged across the 39 TFs) was calculated for our inferred targets, as well as the inferred targets from other methods. We compare the average proportion of over-represented functions between methods in Figure 2. As observed in our knockout validation, our model using expression data only ("Exp Only") does not perform nearly as well as our model with multiple data sources ("All 3" or "Exp+ChIP"). Our model without promoter data ("Exp+ChIP") actually performs better than our model with all three data sources ("All 3"), which is also seen in a subset of the knockout experiments above. Both of these integrated versions of our model have a higher proportion of over-represented functions compared to previous integrated methods, though the performance for all integrated methods are quite similar. The inferred target genes from each of the integrated methods show substantially greater functional over-representation than inferred target genes from the thresholding of a single data source, which again confirms that combining multiple data sources can improve inference. An interesting side note is that the version of our model using expression data alone gives better performance compared with using thresholded expression data. This result suggests that our model for expression data captures additional signal compared to a threshold approach even without integrating additional data sources, though the integrated versions of our model give even better results.

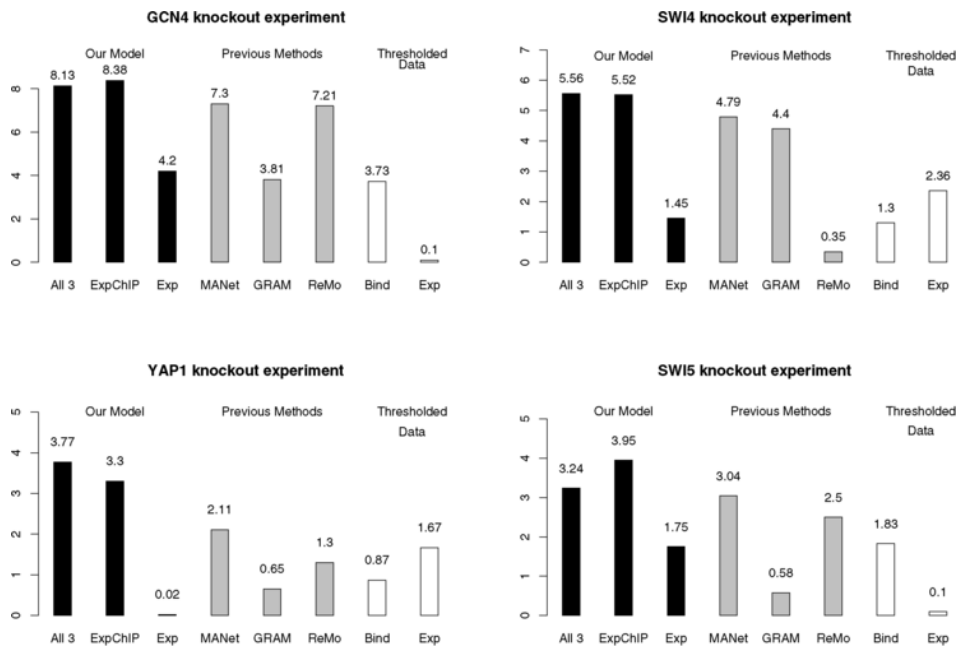

Fig. 1. *T-statistics of knockout response.*

3.2. *Inference for linear model parameters.* Our yeast application involves the simultaneous modeling of multiple transcription factors, which also allows us to infer the partial linear effects $\boldsymbol{\beta}$ of individual transcription factors as well as interaction effects $\boldsymbol{\gamma}$ between pairs of transcription factors. We consider a particular parameter $\beta_j$ or $\gamma_{jk}$ as significant if their 95% posterior interval does not contain zero. It should be noted that we are actually examining the posterior distribution of each $\beta_j$ and $\gamma_{jk}$ parameter conditional on our regulatory indicators $C_{ij} = 1$, since any genes where $C_{ij} = 0$ make no contribution to the conditional distribution of $\beta_j$, as seen in equation (6) of the Appendix. Thus, the parameters $\beta_j$ should be interpreted as the *linear effect of TF j on gene expression if TF j is a regulator of the gene*. Similarly, the parameter $\gamma_{jk}$ should be interpreted as the *interaction effect of TFs j and k on gene expression if TFs j and k are both regulators of the gene*.

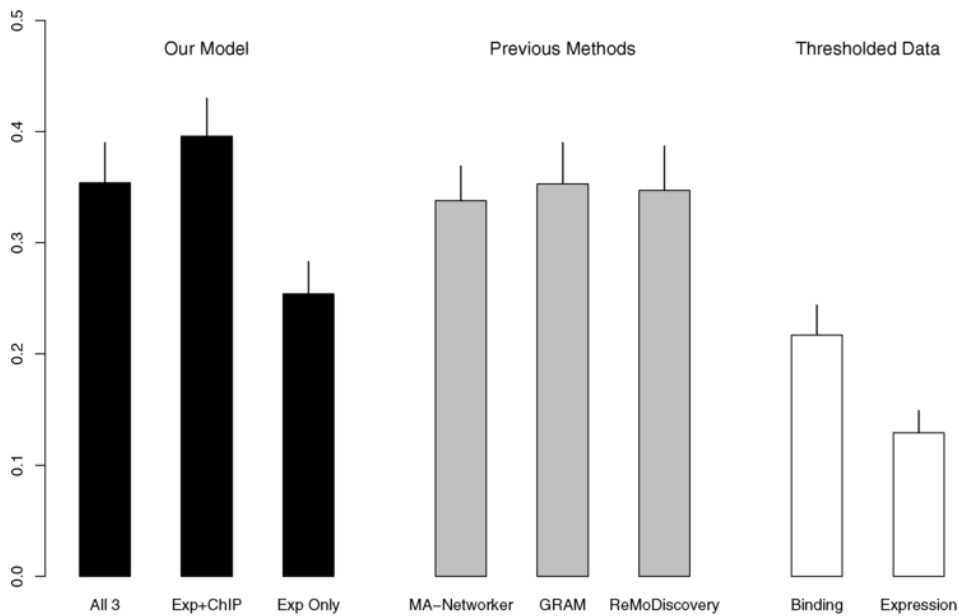

FIG. 2. *Average proportion of genes with over-represented functions. Height of bars represents the average proportion of over-represented functions, while lines represent the standard error of the average.*

Among the linear effects $\boldsymbol{\beta}$, we found sixteen *activators* (significantly positive $\beta_j$'s) and one *repressor* (significantly negative $\beta_j$), which are listed in the supplemental materials. Fourteen of the sixteen activators and the RME1 repressor discovered by our model were previously reported in the SGD database [SGD project (2005)], which gives further evidence that our method is very effective at distinguishing appropriate regulatory relationships. Our model also identified 196 TF pairs which had significant interaction parameters $\gamma_{jk}$. Using our inferred regulation indicators $\boldsymbol{C}$, we imposed an additional restriction that each significant pair of TFs had to also share at least four target genes in common, which resulted in a reduced set of 84 TF pairs. A substantial subset of these 84 TF pairs discovered by our model are also validated by previous biological studies, as outlined in the supplementary materials.

3.3. *Inference about weighting parameters.* A novel component of our proposed methodology was the introduction of a weighting variable $w_j$ which balances the relative quality of the prior distribution based on the ChIP binding data versus the prior distribution based on the promoter sequence data for each TF $j$ individually. Figure 3 gives a boxplot representation of the posterior distributions of the weight variables $w_j$ for all 39 transcription factors.

We see a substantial amount of heterogeneity between each weight variable, which reflects differences in the quality of available data for different transcription factors. We also observe that the posterior distributions for nearly all of these transcription factors are centered around values substantially higher than 0.5, which suggests that the ChIP binding data is being favored as the superior source of prior information for our variable selection indicators $\boldsymbol{C}$. The most extreme example is the RME1 transcription factor, where essentially all of the posterior mass for $w_j$ is greater than 0.8. This general trend matches the common perception of practitioners that a ChIP binding experiment will provide better evidence of regulation than predictions based on sequence data. Other examples include the four transcription factors (indicated by a "K" symbol in Figure 3) for which we have TF knockout data. In Section 3.1 we noticed that the knockout response was generally similar between our full model with all three data sources and a reduced model without the promoter sequence data ($w_j$ set equal to 1 for all $j$). Now we see that this similar performance is expected, considering that the distribution of $w_j$ in the full model is centered quite close to one anyways. It should be noted, however, that this phenomenon is not uniform across all transcription factors. Not all posterior distributions of $w_j$ are pushed toward the boundary value of 1, and in a few cases, such as MSN4, include some posterior mass less than 0.5, which is evidence that the promoter sequence data is also making a contribution to inference.

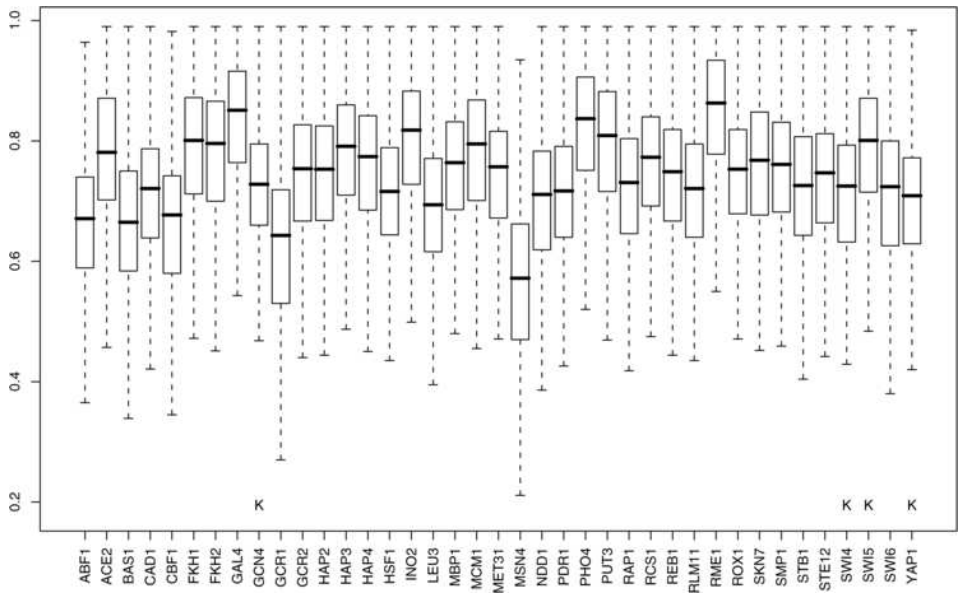

Fig. 3. *Posterior distributions of weight variables.*

## 4. Model sensitivity and consequences.

4.1. *Network sparsity.* Most gene regulatory networks are inherently quite sparse with only a small subset of all genes controlled by any one transcription factor. In terms of our parameterization, this concept translates into an expectation that, for any $j$, only a small number of genes $i$ will have $C_{ij} = 1$. There are a variety of variable selection methods that enforce sparsity on the selection space, such as the lasso [Tibshirani (1996), Efron et al. (2004)]. In the Bayesian variable selection approach, one can also incorporate sparsity by using a small prior probability on the selection indicators, $p(C_{ij} = 1) = \alpha$, where $\alpha$ is small (e.g., $\alpha = 0.01$). In our model, we do not have a constant prior probability $\alpha$ for each selection indicator. Rather, we have specific prior probabilities based our ChIP and sequence data, as in (3). However, as seen in Figure 4, these probabilities $b_{ij}$ and $m_{ij}$ tend to be quite small themselves. To further investigate the sparsity in our model, we repeatedly generated regulatory networks $\boldsymbol{C}$ from our prior distribution (3). From these repeatedly ($m = 10000$) generated networks, we estimated the probability $P(C_{ij} = 1)$ for each gene $i$ and gene $j$, and tabulated the number $N_j$ of inferred target genes for each TF $j$ [genes $i$ with $P(C_{ij} = 1) \geq 0.5$]. This procedure is analogous to the inference from our full model, but only uses the prior probabilities based on ChIP and sequence data. This entire experiment was repeated for different values (ranging from 0.05 to 0.95) of each weight parameter $w_j$, so that we have a range of $N_j$ values for each TF $j$ depending on the different values of $w_j$. Figure 5 shows a boxplot that indicates the range of $N_j$ values over all values $w_j$ for each of our 39 transcription factors. We see that the number of inferred target genes $N_j$ for each TF $j$ is quite small relative to the total number of genes in the network ($\approx 6000$). These results demonstrate that our prior distribution on the network selection indicators $C_{ij}$ is capturing our prior expectations that the Yeast regulatory network should be relatively sparse. We also see substantial differences between TFs in terms of the variability of the number of target genes $N_j$, which is indicative of significance between-TF variability in the response of the inferred target gene set to changing values of the weight $w$. This result provides further motivation for the use of TF-specific weights $w_j$ that balance the ChIP and sequence motif data. Finally, we also included the number of inferred genes $N_j$ from the posterior distribution of our full model (black dots) in Figure 5. We see that our posterior inference differs substantially from our prior inference for most transcription factors, but the overall sparsity of the network is maintained in our posterior distribution.

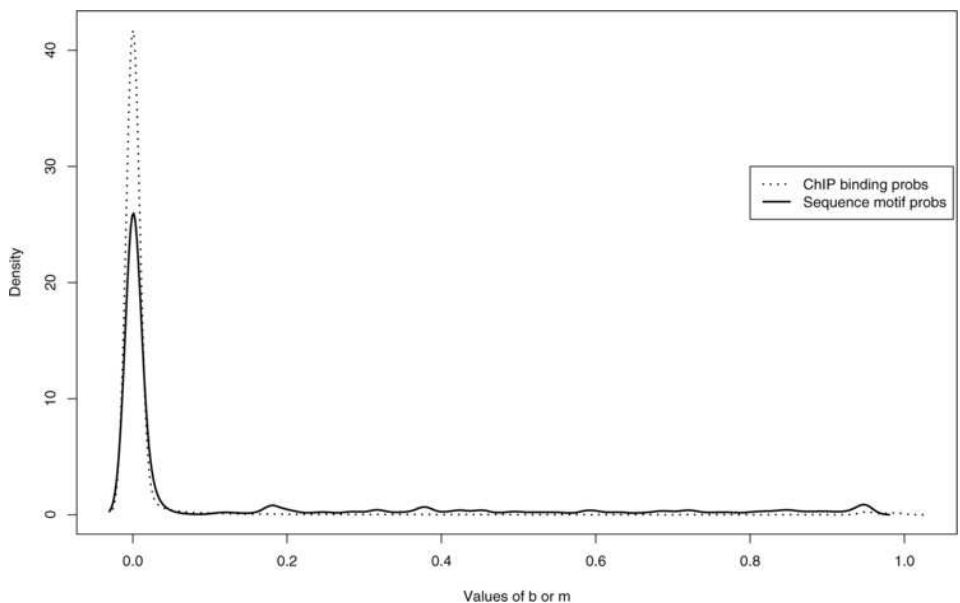

Fig. 4.   *Distribution of $b_{ij}$ and $m_{ij}$ values.*

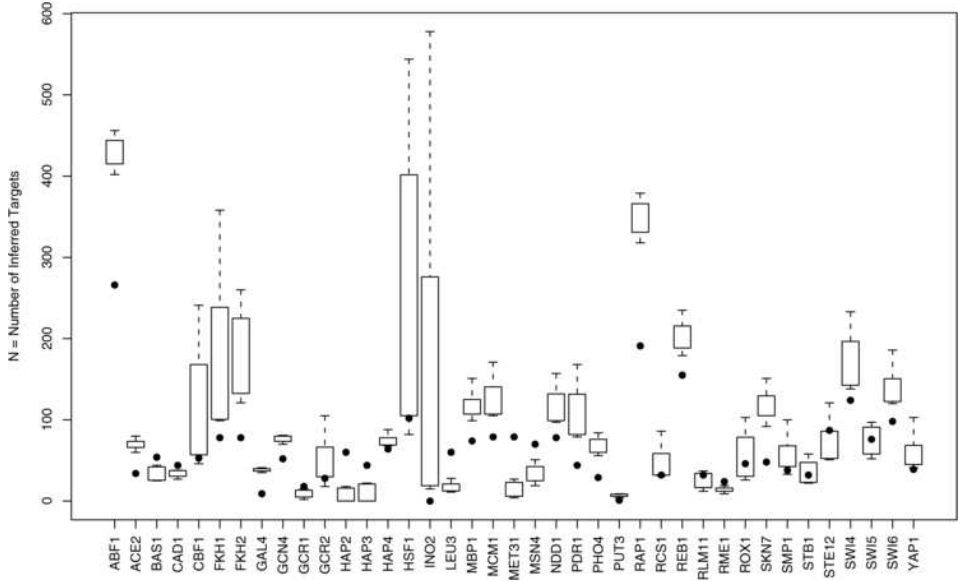

Fig. 5. *Boxplots of distribution of $N_j$ for each TF $j$. Black dots are number of inferred genes from full model.*

4.2. *Sensitivity to prior specification.* Our prior distribution $p(C_{ij}|m_{ij}, b_{ij}, w_j)$ for each regulation indicator $C_{ij}$ is designed to balance the influence of the ChIP binding probability $b_{ij}$ and promoter element probability $m_{ij}$. This balance is achieved in equation (3) by using a weighted geometric mean. We also explored alternative prior specifications that balance our ChIP binding and promoter element data sources. Specifically, we also considered a prior distribution for $C_{ij}$ based on the arithmetic mean of our ChIP binding probability $b_{ij}$ and promoter element probability $m_{ij}$,

$$
\begin{aligned}
(4) \quad & p(C_{ij}|m_{ij}, b_{ij}, w_j) \\
& = [w_j b_{ij} + (1 - w_j) m_{ij}]^{C_{ij}} [1 - (w_j b_{ij} + (1 - w_j) m_{ij})]^{1 - C_{ij}}.
\end{aligned}
$$

For all 39 TFs, we explored differences between the prior probabilities using equation (4) to the prior probabilities using equation (3). Specifically, we examined the difference between the two prior distributions in terms of the number of a priori inferred target genes for each TF [i.e., number of genes $i$ with $p(C_{ij} = 1|m_{ij}, b_{ij}, w_j) \geq 0.5$ for each TF $j$]. For this calculation, we needed to assume a reasonable value for each weight $w_j$, so we used

the posterior mean of each $w_j$ from Section 3.3. Only three TFs (ABF1, RAP1, REB1) showed a substantial difference in the number of *a priori* inferred target genes between the two priors, though these three TFs also had the largest total number of a priori inferred target genes among all 39 TFs. For the remaining TFs, there was very little difference in the number of a priori inferred target genes. We also examined the differences in the a posteriori inferred target genes for the transcription factor HAP4, and found that the list of inferred target genes was quite similar regardless of whether our original prior (3) or the alternative prior (4) were implemented. We evaluated the small differences between the inferred gene lists using our functional categories validation measure (Section 3.1), and found that the inferred target genes using our original prior (3) gave a slightly higher proportion of over-represented functional categories. Given the observed lack of substantial difference between the two prior formulations, and since the specification (4) is a more complicated functional form to implement in our Gibbs sampler, we prefer the use of our original prior specification (3).

Another issue is the potential sensitivity of our posterior inference to the specified prior distributions for the parameters $(\alpha, \beta, \gamma, \sigma^2)$ which appear in the linear model for the expression data (1). The influence of the prior distributions given in Section 2 depends on the values of the hyper-parameters $(\nu, \tau_\alpha^2, \tau_\beta^2, \tau_\gamma^2)$. The nature of the dependence is clear in the conditional distribution formulas in Appendix. The influence of the prior distribution on posterior inference for $\sigma^2$ is very small when $\nu$ is small. The prior distributions for the regression coefficients $\boldsymbol{\alpha}, \boldsymbol{\beta}$ and $\boldsymbol{\gamma}$ can also be made noninfluential by making the prior variance hyperparameters $\tau_\alpha^2, \tau_\beta^2$ and $\tau_\gamma^2$ very large. Our posterior results given in Section 3 are based on values of $\nu = 2$ and $\tau_\alpha^2 = \tau_\beta^2 = \tau_\gamma^2 = 10000$. In many variable selection problems with a large but sparse covariate space, more informative prior distributions on the regression coefficients are used that enforce shrinkage toward values of zero (exclusion of variables). However, that is not necessary in this case, since we have enforced sparsity in our model directly through the selection indicator variables $\boldsymbol{C}$, as detailed in Section 4.1.

**5. Discussion.** We have presented a Bayesian hierarchical model for combining heterogeneous sources of biological data to infer regulatory relationships between genes and transcription factors. Within a variable selection framework, we build upon previous linear models for gene expression data [Bussemaker, Li and Siggia (2001), Tadesse, Vannucci and Lio (2004), Gao, Foat and Bussemaker (2004)] by allowing interactions between transcription factors and incorporating additional information about regulation based on other data sources. The Bayesian paradigm allows us to incorporate these additional data sources in a natural way through the use of prior distributions for our variable selection indicators. This variable selection model also permits genes to belong to multiple regulatory clusters, which allows us to model multiple biological pathways simultaneously. Our full probabilistic model does not rely on any pre-clustering of our data and reduces dependence on arbitrary parameter cutoffs compared to previous methods [e.g., Liao et al. (2003), Yang (2005), Boulesteix and Strimmer (2005)]. When applied to available data in *Saccharomyces cerevisiae* (Yeast), the inferred relationships from our model with multiple data sources were shown to be biologically relevant using external validation measures, with substantially better performance compared with predictions from previous methods (MA-Networker, GRAM and ReMoDiscovery), as well as predictions from thresholding of a single data source. In addition to inferring gene–TF relationships, our model also estimated synergistic and antagonistic interactions between transcription factors, many of which were also validated by previous studies.

The use of informative vs. noninformative prior distributions is a topic of continued discussion within the Bayesian statistical community. Noninformative prior distributions are often used in the context where very little prior information is known, but the researcher still prefers a Bayesian inferential approach for their applied problem. In other cases, prior information is known about the applied problem, in which case the Bayesian paradigm provides a natural way to build this additional information into the probability model. Our current methodology provides a pragmatic compromise of these two approaches: we use informed prior distributions for our primary inferential targets, the regulation indicators $\boldsymbol{C}$, but our model also involves noninformative prior distributions for parameters of secondary interest, such as the coefficients of our linear model for expression data. Our approach of building additional data sources into our model via an informed prior distribution for our regulation indicators $\boldsymbol{C}$ contrasts with most previous Bayesian variable selection research, where criteria are used that assume noninformative prior distributions or avoid prior specification entirely. See [George (2000)] for a review of these noninformative methods and [George and McCulloch (1996)] for a hierarchical Bayesian variable selection model using noninformative prior distributions.

In some previous cases, prior knowledge is incorporated into variable selection, as in the regression model of Garthwaite and Dickey (1996), the logistic regression model of Chen, Ibrahim and Yiannoutsos (1999) and the generalized linear mixed models of Chen et al. (2003). Even more related to our application, Sabatti (2005) used an informed prior distribution based on binding site data to model regulatory networks. However, in their model, only regulation indicators with strong prior evidence are allowed to be nonzero, so that a gene–TF relationship without prior evidence based on sequence data is not permitted regardless of the evidence from gene expression data. Despite relaxing this restriction in our model, our inferred gene regulatory network remains quite sparse, as seen in Section 4.1. The popular "network component analysis" approach [Liao et al. (2003), Yang (2005), Boulesteix and Strimmer (2005)] also assumes that the relationships derived from ChIP binding data alone are known without error. This is a rather restrictive assumption, especially when one considers that ChIP experiments are typically limited to a single condition, but TF binding can vary across different conditions. In contrast, our model allows inferred relationships based on strong gene expression evidence that are not completely evident based on our prior information (ChIP binding or promoter sequence data). Although TF expression is not a perfect proxy for TF activity, we believe it is the best experiment-specific measure of TF activity that our current data resources allow. Our model could certainly be further improved by using a more direct measure of TF activity, such as actual TF protein levels in the cell, but available data on TF protein levels is severly limited at this time.

A fundamental element of our informed prior approach is that we actually have a choice between two prior distributions for our regulation indicators $C$, one informed prior based on ChIP binding data, and another based on promoter sequence data. Since we do not know from application to application (in this case, from transcription factor to transcription factor) which data source is more accurate, we introduce a variable weight that provides a balance between the two prior distributions. This weight variable $w$ is itself assigned a noninformative uniform prior distribution, and we also assign noninformative prior distributions for our linear model parameters. The weighting of different sources of information in a Bayesian model is briefly mentioned by Berry and Hochberg (1999). Ibrahim and Chen (2000) and Chen et al. (2003) introduce the power prior distribution: a weight between their regression model likelihood and a prior distribution based on historical data. In contrast, our weight variables are used as a balance between two "competing" prior distributions, which means that the estimated posterior distribution of each weight variable can shed substantial insight into the relative quality of our two sources of prior data. In fact, our weighted prior distribution can be interpreted as the combination of our two sources of prior information that best matches the likelihood distribution based on expression data. The results from our Yeast application indicate that our variable weight methodology achieves an appropriate balance between our two sources of prior information. Our results confirm the commonly-held belief that promoter sequence data is generally much less reliable than the ChIP binding data, although promoter sequence data can be useful in some cases.

## APPENDIX: GIBBS SAMPLING IMPLEMENTATION

The posterior distribution of our unknown parameters is proportional to the product of our model likelihood and our assumed prior distributions,

$$p(\boldsymbol{C}, \boldsymbol{w}, \boldsymbol{\Theta} | \boldsymbol{g}, \boldsymbol{f}, \boldsymbol{m}, \boldsymbol{b})$$

$$\propto p(\boldsymbol{g} | \boldsymbol{f}, \boldsymbol{C}, \boldsymbol{\Theta}) \cdot p(\boldsymbol{C} | \boldsymbol{m}, \boldsymbol{b}, \boldsymbol{w}) \cdot p(\boldsymbol{w}) \cdot p(\boldsymbol{\alpha}) \cdot p(\boldsymbol{\beta}) \cdot p(\sigma^2)$$

$$= \prod_{i=1}^{N} \prod_{t=1}^{T} (2\pi\sigma^2)^{-1/2} \exp\left[\frac{-1}{2\sigma^2}\left(g_{it} - \alpha_i - \sum_{j=1}^{J} \beta_j X_{ijt}\right)^2\right]$$

$$\times \prod_{i=1}^{N} \prod_{j=1}^{J} [b_{ij}^{C_{ij}}(1-b_{ij})^{1-C_{ij}}]^{w_j} \cdot [m_{ij}^{C_{ij}}(1-m_{ij})^{1-C_{ij}}]^{1-w_j}$$

$$\times \prod_{i=1}^{N} \frac{1}{\tau_\alpha} \exp\left[\frac{-1}{2\tau_\alpha^2}\alpha_i^2\right] \cdot \prod_{j=1}^{J} \frac{1}{\tau_\beta} \exp\left[\frac{-1}{2\tau_\beta^2}\beta_j^2\right] \cdot \frac{1}{\sigma} \exp\left[\frac{-1}{2\sigma^2}\right],$$

where $X_{ijt} = C_{ij}f_{jt}$. We use the following Gibbs sampling [Geman (1984)] steps to estimate the joint posterior distribution of all unknown parameters.

*Step* 1. *Sampling linear model parameters* $\boldsymbol{\Theta}$. The regulation matrix $\boldsymbol{C}$ is assumed known during this step, so we do not need to use our prior data $\boldsymbol{b}, \boldsymbol{m}$ or the current values of $\boldsymbol{w}$. We use $\boldsymbol{C}$ to construct the variables $\boldsymbol{X}$, where $X_{ijt} = C_{ij}f_{jt}$. The linear model parameters $\boldsymbol{\Theta}$ are then separately estimated by the following iterative strategy. Note that in the steps below, we have combined our interaction coefficients $\gamma_{jk}$ and linear coefficients $\beta_j$ into a single set of parameters $\boldsymbol{\beta}$. Since each intercept $\alpha_i$ is independent from the other $\alpha$'s, they can be separately sampled,

$$p(\alpha_i | \boldsymbol{\beta}, \sigma^2, \boldsymbol{g}, \boldsymbol{X}) \propto \exp\left[\frac{-1}{2\sigma^2}\sum_{t=1}^{T}\left(g_{it} - \alpha_i - \sum_{j=1}^{J}\beta_j X_{ijt}\right)^2\right] \cdot \exp\left[\frac{-1}{2\tau_\alpha^2}\alpha_i^2\right]$$

$$= \exp\left[\frac{-1}{2\nu_\alpha}\left(\alpha_i - \frac{\nu_\alpha}{\sigma^2}\cdot\sum_{t=1}^{T}Y_t\right)^2\right],$$

where $Y_t = g_{it} - \sum_{j=1}^{J}\beta_j X_{ijt}$ and $\nu_\alpha = (T/\sigma^2 + 1/\tau_\alpha^2)^{-1}$. This distribution implies that

$$(5) \qquad\qquad \alpha_i \sim \text{Normal}\left(\frac{\nu_\alpha}{\sigma^2}\cdot\sum_{t=1}^{T}Y_t, \nu_\alpha\right).$$

We can make our prior distribution for each $\alpha_i$ to be noninformative by making $\tau_\alpha$ very large (in this study, 10000) relative to the contribution of the likelihood to the variance ($\sigma^2/T$).

Our slope and interaction coefficients $\beta_j$'s are not independent from each other, and so must be iteratively sampled:

$$p(\beta_j | \boldsymbol{\alpha}, \sigma^2, \boldsymbol{g}, \boldsymbol{X})$$

$$\propto \exp\left[\frac{-1}{2\sigma^2} \sum_{i=1}^{N} \sum_{t=1}^{T} \left(g_{it} - \alpha_i - \sum_{j'=1}^{J} \beta_{j'} X_{ij't}\right)^2\right] \cdot \exp\left[\frac{-1}{2\tau_\beta^2} \beta_j^2\right]$$

$$= \exp\left[\frac{-1}{2\sigma^2} \sum_{i=1}^{N} \sum_{t=1}^{T} (V_{it} - \beta_j X_{ijt})^2\right] \cdot \exp\left[\frac{-1}{2\tau_\beta^2} \beta_j^2\right],$$

where $V_{it} = g_{it} - \alpha_i - \sum_{j' \neq j} \beta_{j'} X_{ij't}$, which reduces further to

$$p(\beta_j | \boldsymbol{\alpha}, \sigma^2, \boldsymbol{g}, \boldsymbol{X}) \propto \exp\left[\frac{-1}{2\nu_\beta} \left(\beta_j - \frac{\nu_\beta}{\sigma^2} \cdot T_{VX}\right)^2\right],$$

where $\nu_\beta = (T_{XX}/\sigma^2 + 1/\tau_\beta^2)^{-1}$, $T_{XX} = \sum_{i=1}^{N} \sum_{t=1}^{T} X_{ijt}^2$ and $T_{VX} = \sum_{i=1}^{N} \sum_{t=1}^{T} V_{it} X_{ijt}$.

This distribution implies that

$$(6) \qquad\qquad \beta_j \sim \text{Normal}\left(\frac{\nu_\beta}{\sigma^2} \cdot T_{VX}, \nu_\beta\right).$$

We can make our prior distribution for each $\beta_j$ to be noninformative by making $\tau_\beta$ very large (in this study, 10000) relative to the contribution of the likelihood to the variance $(\sigma^2/T_{XX})$.

For the residual variance $\sigma^2$, we use a $\chi_\nu^2$ prior distribution for $\sigma$ with hyper-parameter $\nu = 2$, which results in the following conditional distribution:

$$p(\sigma^2 | \boldsymbol{\alpha}, \boldsymbol{\beta}, \boldsymbol{g}, \boldsymbol{X})$$

$$\propto (\sigma^2)^{-(TN/2+2)} \cdot \exp\left[\frac{-1}{2\sigma^2} \sum_{i=1}^{N} \sum_{t=1}^{T} \left(g_{it} - \alpha_i - \sum_{j=1}^{J} \beta_j X_{ijt}\right)^2\right]$$

$$= (\sigma^2)^{-((TN+2)/2+1)} \cdot \exp\left[\frac{-1}{2\sigma^2} (V_\sigma + 1)^2\right],$$

where $V_\sigma = \sum_{i=1}^{N} \sum_{t=1}^{T} (g_{it} - \alpha_i - \sum_{j=1}^{J} \beta_j X_{ijt})^2$. We see that the influence of this prior is very small on the posterior distribution for $\sigma^2$, which is a scaled-inverse $\chi^2$ distribution with degrees of freedom parameter $TN + 2$ and scale parameter $s^2 = (V_\alpha + 1)/(TN + 2)$.

*Step* 2. *Sampling regulation indicators* $\boldsymbol{C}$. We are assuming that both our linear model parameters $\boldsymbol{\Theta}$ and our weights $\boldsymbol{w}$ are known for this step of the algorithm. When estimating a new value for each $C_{ij}$, we also can condition on $\boldsymbol{C}'$, which is all the other $C_{i'j'}$ values in $\boldsymbol{C}$ ($i' \neq i$ and $j' \neq j$). This gives us the following conditional distribution for $C_{ij}$:

$$
\begin{aligned}
(7) \quad & p(C_{ij}|\boldsymbol{\Theta}, \boldsymbol{w}, \boldsymbol{C}', \boldsymbol{g}, \boldsymbol{f}, \boldsymbol{b}, \boldsymbol{m}) \\
& \propto \exp\left[\frac{-1}{2\sigma^2}\sum_{t=1}^{T}\left(g_{it} - \alpha_i - \sum_{j=1}^{J}\beta_j C_{ij} f_{jt}\right)^2\right] \\
& \quad \times [b_{ij}^{C_{ij}}(1-b_{ij})^{1-C_{ij}}]^{w_j} \cdot [m_{ij}^{C_{ij}}(1-m_{ij})^{1-C_{ij}}]^{1-w_j}.
\end{aligned}
$$

Let $Z_1$ be the value of equation (7) when $C_{ij} = 1$ and $Z_0$ be the value of equation (7) when $C_{ij} = 0$. We sample a new value of $C_{ij}$ equal to 1 or 0 with probabilities proportional to $Z_1$ or $Z_0$ respectively.

*Step* 3. *Sampling prior weights* $\boldsymbol{w}$. We are assuming that the regulation matrix $\boldsymbol{C}$ is known for this step of the algorithm, so we do not need to use any of the expression data, $\boldsymbol{g}$ or linear model parameters $\boldsymbol{\Theta}$ for this step. For each TF $j$, we need to sample a new weight $w_j$ based on the following distribution:

$$
\begin{aligned}
(8) \quad & p(w_j|\boldsymbol{C}, \boldsymbol{b}, \boldsymbol{m}) \propto A(w_j)^{-1} \cdot \left[\prod_{i=1}^{n}(b_{ij})^{C_{ij}}(1-b_{ij})^{1-C_{ij}}\right]^{w_j} \\
& \quad \times \left[\prod_{i=1}^{n}(m_{ij})^{C_{ij}}(1-m_{ij})^{1-C_{ij}}\right]^{1-w_j}.
\end{aligned}
$$

The normalizing constant $A(w_j)$ present in (8) comes from the integration of $p(w_j, \boldsymbol{C}_j|\boldsymbol{b}, \boldsymbol{m})$ over all configurations of $\boldsymbol{C}_j$:

$$
(9) \quad A(w_j) = \prod_{i=1}^{n}[(b_{ij})^{w_j}(m_{ij})^{1-w_j} + (1-b_{ij})^{w_j}(1-m_{ij})^{1-w_j}].
$$

We sample a new value $w_j$ via *grid sampling*: we evaluate (8) over a fine grid of points in the unit interval, and sample one of these points with probability proportional to (8). Multiple chains of our Gibbs sampling algorithm were run from different starting points until we were confident that the chains had converged to the same range of values. Details of our convergence diagnostics are given in the supplemental materials.

**Acknowledgments.** We would like to thank Jonathan Schug and Edward I. George for helpful discussions.

S. T. JENSEN
DEPARTMENT OF STATISTICS
THE WHARTON SCHOOL
UNIVERSITY OF PENNSYLVANIA
PHILADELPHIA, PENNSYLVANIA 19104
USA
E-MAIL: stjensen@wharton.upenn.edu

G. CHEN
DEPARTMENT OF BIOENGINEERING
CENTER FOR BIOINFORMATICS
UNIVERSITY OF PENNSYLVANIA
PHILADELPHIA, PENNSYLVANIA 19104
USA
E-MAIL: ggchen@pcbi.upenn.edu

C. J. STOECKERT, JR.
DEPARTMENT OF GENETICS
CENTER FOR BIOINFORMATICS
UNIVERSITY OF PENNSYLVANIA
PHILADELPHIA, PENNSYLVANIA 19104
USA
E-MAIL: stoeckrt@pcbi.upenn.edu